\documentclass[12pt]{article}
\usepackage{amsmath,amsfonts,amssymb,amsthm}
\usepackage{mathrsfs}
\usepackage{latexsym}
\usepackage[english]{babel}

\usepackage{graphicx}
\usepackage[usenames,dvipsnames]{xcolor}

\renewenvironment{proof}[1][\proofname]{{\bfseries #1.} }{\qed}

\def\Cov{{\rm Cov\,}}

\newcommand{\field}[1]{\mathbb{#1}}

\newcommand{\R}{\field{R}}

\newcommand{\Z}{\field{Z}}

\newcommand{\Var}{{\rm Var}}

\newcommand{\e}{{\rm e}}

\newcommand{\eps}{\varepsilon}

\def\authors#1{{ \begin{center} #1 \vspace{0pt} \end{center} } \smallskip}
\def\institution#1{{\sl \begin{center} #1 \vspace{0pt} \end{center} } }

\def\title#1{{\huge\bf  \begin{center} #1 \vspace{0pt} \end{center}  } \smallskip}

\def\E{{\mathbb{ E}}}

\def\P{{\mathbb{P}}}

\def\paref#1{(\ref{#1})}

\newtheorem{theorem}{Theorem}[section]

\newtheorem{definition}[theorem]{Definition}

\newtheorem{remark}[theorem]{Remark}

\begin{document}

\date{Jan 2018}

\title{\sc {Random Nodal Lengths and\\ Wiener Chaos}}
\authors{\large Maurizia Rossi
}
\institution{
MAP5-UMR CNRS 8145, Universit\'e Paris Descartes, France \\
E-mail address: \texttt{maurizia.rossi@parisdescartes.fr}
}

\begin{abstract}

In this survey we collect some of the recent results on the ``nodal geometry" of random eigenfunctions on Riemannian surfaces. We focus on the asymptotic behavior, for high energy levels, of the nodal length of Gaussian Laplace eigenfunctions on the torus (arithmetic random waves) and on the sphere (random spherical harmonics).  We give some insight on both Berry's cancellation phenomenon and the nature of nodal length second order fluctuations (non-Gaussian on the torus and Gaussian on the sphere) in terms of chaotic components. 
Finally we consider the general case of monochromatic random waves, i.e. Gaussian random linear combination of eigenfunctions of the Laplacian on a compact Riemannian surface with frequencies from a short interval, whose scaling limit is Berry's Random Wave Model. For the latter we present some recent results on the asymptotic distribution of its nodal length in the high energy limit (equivalently, for growing domains). 

\smallskip


\smallskip

\end{abstract}

\section{Introduction}

\subsection{Background and motivations}\label{BG}

\begin{figure}[htbp]
\centering
\includegraphics[width=7cm]{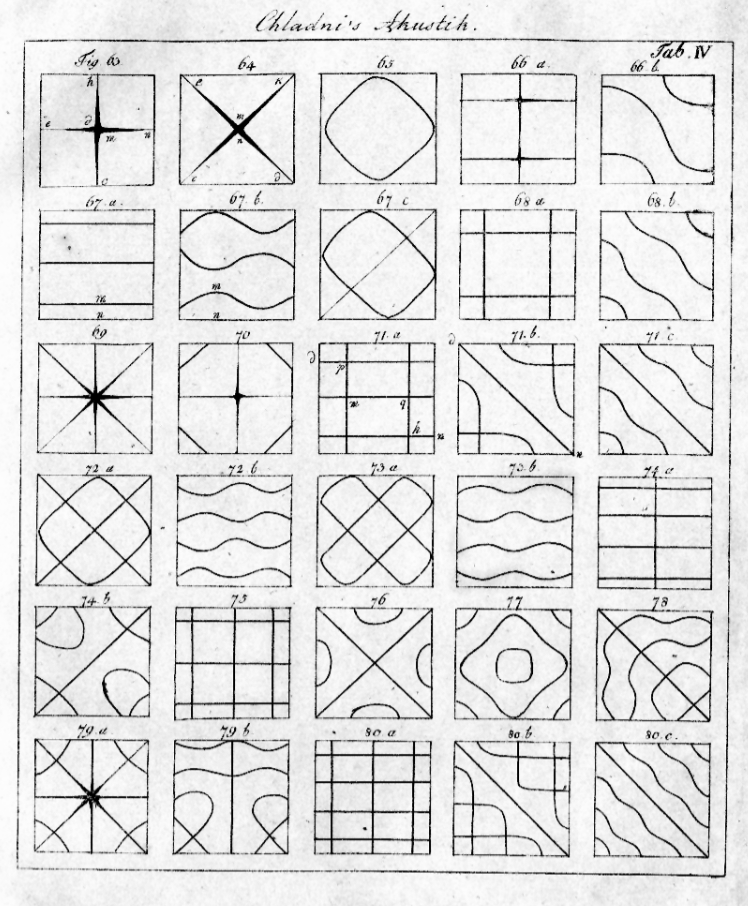}
\caption{Chladni's figures\label{fig1}}
\end{figure}

At the end of the 18th century, the musician and physicist
Chladni noticed that  sounds of different
pitch could be made by exciting a metal plate with the bow of a violin, depending on where the bow touched the plate. The latter was
fixed only in the center, and when there was some sand on it, for
each pitch a curious pattern appeared (see Figure \ref{fig1}).
About 60 years later, Kirchhoff, inspired also by previous contributions of Germain, Lagrange and Poisson, showed that Chladni figures on a plate correspond to the zeros of eigenpairs
(eigenvalues and corresponding eigenfunctions) of the biharmonic operator with free
boundary conditions. Nowadays nodal patterns arise in several areas, from  the musical instruments industry to the study of natural phemomena such as earthquakes. 
See \cite{Chl, GK12, Wig12} and the references therein for more details. 

Now let $(\mathcal M,g)$ denote a compact smooth Riemannian manifold of dimension $2$ and $f:\mathcal M\to \mathbb R$ a function on $\mathcal M$. The set $f^{-1}(0) := \lbrace x\in \mathcal M : f(x) = 0\rbrace$ is usually called \emph{nodal set} of $f$. 
We are interested in the case where $f$ is an eigenfunction of the  Laplace-Beltrami operator 
 $\Delta$ on the manifold. 
There exists an orthonormal basis $\lbrace f_j \rbrace_{j=0}^{+\infty}$ of $L^2(\mathcal M)$ consisting of eigenfunctions for $-\Delta$ whose corresponding sequence of  eigenvalues $\lbrace \lambda_j^2 \rbrace_{j=0}^{+\infty}$
\begin{equation}\label{eq}
-\Delta f_j = \lambda^2_j f_j
\end{equation} 
is non-decreasing (in particular $0=\lambda_0 < \lambda_1 \le \lambda_2 \le \dots)$. 

The nodal set of $f_j$  is, generically, a smooth curve (called \emph{nodal line}) whose components are homeomorphic to the circle. We are interested in the geometry of $f_j^{-1}(0)$ for high energy levels, i.e. as $\lambda_j\to +\infty$.  

Yau's conjecture \cite{Yau82, Yau93} concerns, in particular\footnote{Yau stated its conjecture in any dimension}, the length of nodal lines $f_j^{-1}(0)$:  there exist two positive constant $c_\mathcal M, C_\mathcal M$ such that 
\begin{equation}\label{yau_conj}
c_{\mathcal M} \lambda_j \le \text{length}(f_j^{-1}(0)) \le C_{\mathcal M} \lambda_j,
\end{equation}
for every $j\ge 1$. Yau's conjecture was proved by Donnelly and Fefferman in \cite{DF88} for real analytic manifolds, and the lower bound
was established by Logunov and Malinnikova \cite{Log16a, Log16b, LM16} for the general case. 

In \cite{Ber77} Berry conjectured that the local behavior of high energy eigenfunctions $f_j$ in \paref{eq} for generic \emph{chaotic} surfaces is universal. He meant that it is comparable to the behavior of the centred Gaussian field $B_{\lambda_j} = \lbrace B_{\lambda_j}(x)\rbrace_{x \in \R^2}$ on the Euclidean plane whose covariance kernel is, for $x,y\in \R^2$, 
\begin{equation}\label{cov_berry}
\Cov(B_{\lambda_j}(x), B_{\lambda_j}(y))=J_0(\lambda_j \| x - y\|),
\end{equation}
$J_0$ being the Bessel function of order zero \cite[\S 1.71]{szego} and $\| \cdot \|$ denoting the Euclidean norm (the latter model is known as Berry's Random Wave Model). In particular, nodal lines of $B_{\lambda_j}$ should model nodal lines of deterministic eigenfunctions $f_j$ for large eigenvalues, allowing the investigation of associated local quantities like the (nodal) length. 
See \cite{Wig12} and the references therein for more details.

\subsection{Short plan of the survey}

In the case of surfaces $(\mathcal M, g)$  with spectral degeneracies, like the two dimensional standard flat torus or the unit round sphere, it is possible to construct a random model by endowing each eigenspace with a Gaussian measure. The asymptotic behavior, in the high energy limit, of the corresponding nodal length is the subject of \S \ref{Snodal_results}. In particular, two comments are in order: the asymptotic variance is smaller than expected in both cases (\emph{Berry's cancellation phenomenon}), and the nature of second order fluctuations is Gaussian for the spherical case but not on the torus. In \S \ref{Sgreen} we give some insight on both these phenomena in terms of \emph{chaotic projections}, notion introduced in \S \ref{Swiener}. 

For a generic Riemannian metric on a smooth compact surface, the eigenvalues  of the Laplace-Beltrami operator are simple and one cannot define a meaningful random model as for the toral or the spherical case, that is, a single frequency random model. In view of this, in \S \ref{Sfurther} 
we work with the so-called \emph{monochromatic random waves} i.e., Gaussian random linear combination of eigenfunctions  with frequencies from a short interval. Their scaling limit is Berry's RWM \paref{cov_berry}, and for the latter we present some recent results on the asymptotic behavior for its nodal length, which turns out to be Gaussian as for the spherical case. 

\subsection*{Acknowledgements} The author would like to thank the organizers of the conference \emph{Probabilistic methods in spectral geometry and PDE}  and the CRM Montr\'eal for such a wonderful meeting. 
Some of the results presented in this short survey have been obtained in collaboration with Domenico Marinucci, Ivan Nourdin, Giovanni Peccati and Igor Wigman to whom the author is grateful; in particular, the author would like to thank Domenico Marinucci for the insight in \S \ref{Sgreen} and Igor Wigman for fruitful discussions we had in Montr\'eal. Moreover, the author thanks Valentina Cammarota for a number of useful discussions, and an anonymous referee for his/her valuable comments. 
The author was financially supported by the grant F1R-MTH-PUL-15STAR
(STARS) at the University of Luxembourg and the ERC grant no. 277742 \emph{Pascal}.
Her research is currently supported  by the 
Fondation Sciences Math\'ematiques de Paris and the project ANR-17-CE40-0008 \emph{Unirandom}.

\section{Random nodal lengths}\label{Snodal_results}

\subsection{The toral case}\label{SresultsT}

\subsubsection{Arithmetic random waves} 
The Laplacian eigenvalues for the two-dimensional standard flat torus $\mathbb T^2:=\R^2/\Z^2$ are of the form $4\pi^2n$, where $n$  is an integer expressible as a sum of two squares
$
(n\in S := \lbrace n = a^2 + b^2 : a,b\in \Z \rbrace)
$.
For $n\in S$, we define $\Lambda_n$ to be the set of frequencies 
\begin{equation}\label{Lambda}
\Lambda_n := \lbrace \xi\in \mathbb Z^2 : \| \xi\|= \sqrt n\rbrace 
\end{equation}
and we denote by $\mathcal N_n$ its cardinality ($\mathcal N_n$ is the multiplicity of  eigenvalue $4\pi^2n$). The set $\Lambda_n$ in \paref{Lambda} induces a probability measure $\mu_n$  on the unit circle $\mathbb S^1\subset \R^2$
$$
\mu_n := \frac{1}{\mathcal N_n} \sum_{\xi\in \Lambda_n} \delta_{\xi/\sqrt n},
$$ 
$\delta_z$ denoting the Dirac mass at $z\in \R^2$. For more details see \cite{KKW}.

For $n\in S$, the arithmetic random wave $T_n$ (of order $n$) is the Gaussian random eigenfunction on the torus defined as follows \cite{RW}:
\begin{equation}\label{defARW}
T_n(x) := \frac{1}{\sqrt{\mathcal N_n}} \sum_{\xi\in \Lambda_n}a_\xi \e^{i2\pi \langle \xi, x\rangle},\qquad x\in \mathbb T^2,
\end{equation}
where $\lbrace a_\xi \rbrace_{\xi\in \Lambda_n}$ is a family of i.d. standard complex Gaussian random variables\footnote{defined on some probability space $(\Omega, \mathcal F, \mathbb P)$}, and independent except for the relation $\overline{a_\xi} = a_{-\xi}$ that ensures $T_n$ to be real. Recall that $\lbrace \e^{i2\pi \langle \xi, \cdot \rangle} \rbrace_{\xi\in \Lambda_n}$ is an orthonormal basis for the eigenspace related to the eigenvalue $4\pi^2 n$. 
Equivalently, it can be defined as the centered Gaussian field on $\mathbb T^2$ whose covariance kernel  is, for $x,y\in \mathbb T^2$,
\begin{equation}\label{cov_int}
\Cov(T_n(x), T_n(y)) = \frac{1}{\mathcal N_n} \sum_{\xi\in \Lambda_n} \e^{i 2\pi\langle \xi, x-y\rangle }. 
\end{equation}
There exists a density one subsequence $\lbrace n_j\rbrace_j\subset S$ of energy levels such that, as $j\to +\infty$, 
$$
\mu_{n_j} \Rightarrow d\theta/2\pi,
$$
$d\theta$ denoting the uniform measure on $\mathbb S^1$ (see \cite{FKW06}). 
From \paref{cov_int} we have, for $x,y\in \mathbb T^2$, as $j\to +\infty$, 
\begin{equation*}
\Cov(T_{n_j}(x/2\pi \sqrt{n_j}), T_{n_j}(y/2\pi\sqrt{n_j})) = \int_{\mathbb S^1} \e^{i \langle \theta, x-y\rangle} d\mu_{n_j}(\theta) \to J_0(\|x-y\|),
\end{equation*}
$J_0$ still denoting the Bessel function of order zero (cf. \paref{cov_berry}).
 There exist other weak-$\star$ partial limits of the sequence $\lbrace \mu_n\rbrace_{n\in S}$, partially classified in \cite{KW}. 

\subsubsection{Nodal lengths: some recent results}

The nodal set $T_n^{-1}(0)$ is a smooth curve a.s.; let us set $\mathcal L_n := \text{length}(T_n^{-1}(0))$. By means of Kac-Rice formulas, the expected nodal length was computed by Rudnick and Wigman in \cite{RW}
\begin{equation}\label{mean_torus}
\E[\mathcal L_n] = \frac{1}{2\sqrt 2} \sqrt{4\pi^2n},
\end{equation}
and in \cite{KKW} the exact asymptotic variance was found: as $\mathcal N_n\to +\infty$, 
\begin{equation}\label{varTorus}
\Var(\mathcal L_n) \sim \frac{1+\widehat{\mu_n}(4)^2}{512} \, \frac{4\pi^2n}{\mathcal N_n^2},
\end{equation}
$\widehat{\mu_n}(4)$ denoting the fourth Fourier coefficients of $\mu_n$. In order to have an asymptotic law for the variance, it suffices to choose a subsequence $\lbrace n_j\rbrace_j\subset S$  such that as $j\to +\infty$ we have (a) $\mathcal N_{n_j}\to +\infty$ and (b) $|\widehat{\mu_{n_j}}(4)|\to \eta$,  for some $\eta\in [0,1]$. Note that for each $\eta\in [0,1]$, there exists a subsequence $\lbrace n_j\rbrace$ such that both (a) and (b) hold (see \cite{KKW, KW}). 
The second order fluctuations of $\mathcal L_n$ were investigated and fully resolved in \cite{MPRW}, as follows.  
\begin{theorem}\label{thTorus}
For $\lbrace n_j\rbrace_j\subset S$  such that, as $j\to +\infty$, $\mathcal N_{n_j}\to +\infty$ and $|\widehat{\mu_{n_j}}(4)|\to \eta$  for some $\eta\in [0,1]$, we have 
\begin{equation}\label{NCLT}
\frac{\mathcal L_{n_j} - \E[\mathcal L_{n_j}]}{\sqrt{\Var(\mathcal L_{n_j})}}\mathop{\to}^d \frac{1}{2\sqrt{1+\eta^2}} (2 - (1-\eta) Z_1^2 - (1+\eta)Z_2^2),
\end{equation}
where $Z_1$ and $Z_2$ are i.i.d. standard Gaussian random variables. 
\end{theorem}
A quantitative version (in Wasserstein distance) of \paref{NCLT} is given in \cite{PR}, where an alternative proof for \paref{varTorus} is given by means of chaotic expansions.  In the recent paper \cite{BMW17} the authors study the nodal length of arithmetic random waves
restricted to decreasing domains (shrinking radius-$s$ balls)  all the way down to Planck scale (i.e., for $s> n^{-1/2+\varepsilon}$, $\varepsilon >0$). Remarkably, they prove that the latter is asymptotically fully correlated with $\mathcal L_n$. In \cite{DNPR} the interesection number of nodal lines corresponding to two independent arithmetic random waves with the same eigenvalue is studied, whereas 
the papers \cite{R-W, RoW} investigate the intersection number of the nodal lines $T_n^{-1}(0)$ and a fixed reference curve with nowhere zero curvature. The case of a straight line segment is studied in \cite{Maf16a} by Maffucci.  There are results also in the three dimensional setting \cite{BM17, Cam17, Maf16b, RWY}.

\subsection{The spherical case}\label{SresultsS}

\subsubsection{Random spherical harmonics}
  The Laplacian eigenvalues on the two-dimensional unit round sphere $\mathbb S^2$ are of the form $\ell(\ell+1)$, where $\ell\in \mathbb N$, and the multiplicity of the $\ell$-th eigenvalue is $2\ell+1$. The $\ell$-th random spherical harmonic on $\mathbb S^2$ is the Gaussian random eigenfunction $T_\ell$ (abusing notation) defined as follows (see e.g. \cite{Wig}):
\begin{equation}\label{defSphere}
T_\ell(x) := \sqrt{\frac{4\pi}{2\ell+1}} \sum_{m=-\ell}^\ell a_{\ell,m} Y_{\ell,m}(x),\qquad x\in \mathbb S^2,
\end{equation}
where $\lbrace a_{\ell,m} \rbrace_{m=-\ell, \dots, \ell}$ is a family of i.d. standard complex Gaussian random variables, and independent except for the relation $\overline{a_{\ell,m}} = (-1)^\ell a_{\ell,-m}$ that ensures $T_\ell$ to be real. The family $\lbrace Y_{\ell,m}\rbrace_{m=-\ell,\dots , \ell}$ of spherical harmonics \cite[\S 3.4]{MP11} is an orthonormal basis for the eigenspace related to the eigenvalue $\ell(\ell+1)$. 
Equivalently, it can be defined as the centered Gaussian field on the sphere whose covariance kernel is, for $x,y\in \mathbb S^2$,
$$
\Cov(T_\ell(x)), T_\ell(y)) = P_\ell(\cos d(x,y)),
$$
where $P_\ell$ denotes the $\ell$-th Legendre polynomial \cite[\S 13.1.2]{MP11} and $d(x,y)$ the geodesic distance between the two points $x$ and $y$. 
Hilb's asymptotic formula \cite[Theorem 8.21.12]{szego} states that, \emph{uniformly} for  $\theta\in [0, \pi - \varepsilon]$ ($\varepsilon >0$), as $\ell\to +\infty$, 
\begin{equation}\label{Hilb}
P_\ell(\cos \theta) \sim \sqrt{\frac{\theta}{\sin \theta}} J_0((\ell+1/2)\theta),
\end{equation}
cf. \paref{cov_berry}.

\subsubsection{Nodal lengths: some recent results}
The nodal set $T_\ell^{-1}(0) := \lbrace x\in \mathbb S^2 : T_\ell(x) = 0\rbrace$ is a smooth curve a.s.
The mean of the nodal length (abusing notation) $\mathcal L_\ell := \text{length}(T_\ell^{-1}(0))$ was computed in \cite{Ber85} 
\begin{equation}\label{mean_sphere}
\E[\mathcal L_\ell] = 4\pi\cdot  \frac{1}{2\sqrt 2} \sqrt{\ell(\ell+1)};
\end{equation}
note that the coefficient $1/2\sqrt 2$ in \paref{mean_sphere} is universal (cf. \paref{mean_torus}). 
The asymptotic behaviour of the variance was given in \cite{Wig}: as $\ell\to +\infty$, 
\begin{equation}\label{varSphere}
\Var(\mathcal L_\ell) \sim \frac{1}{32}\log \ell.
\end{equation}
The second order fluctuations of $\mathcal L_\ell$ have been investigated and fully resolved in \cite{MRW17}:
\begin{theorem}\label{thSphere}
As $\ell\to +\infty$, 
$$
\frac{\mathcal L_\ell - \E[\mathcal L_\ell]}{\sqrt{\Var(\mathcal L_\ell)}}\mathop{\to}^d Z,
$$
where $Z$ is a standard Gaussian random variable. 
\end{theorem}

\subsection{Some comments}\label{Sremarks}

\begin{remark}\label{rem1}\rm

(i) The mean of the nodal length in both the toral \paref{mean_torus} and the spherical case \paref{mean_sphere} is of the form $1/2\sqrt 2$ (a universal constant) times the square root of the corresponding eigenvalue (cf \paref{yau_conj}) times the area of the surface.

\noindent (ii) The order of magnitude of the asymptotic variance of the nodal length in both the toral  \paref{varTorus} and the spherical setting \paref{varSphere} is smaller than expected. Indeed, Kac-Rice formula \cite[Ch. 12]{AT} suggests that its variance is asymptotically proportional to the corresponding eigenvalue times the second moment of the covariance kernel which is, up to a constant factor, the inverse of the eigenvalue's multiplicity ($n/\mathcal N_n$ on $\mathbb T^2$ and $\ell$ for $\mathbb S^2$).  This fact is known as Berry's \emph{cancellation phenomenon}, indeed it was predicted by Berry in \cite{berry2002} and proved by Wigman in \cite{Wig} on the sphere and in \cite{KKW} on the torus. 

\noindent (iii) Theorem \ref{thTorus} reveals in particular the non-Gaussian asymptotic nature of the toral nodal length, in contrast with the Gaussian fluctuations for the spherical case (Theorem \ref{thSphere}). 
\end{remark}
We will give some insights for this remark in \S \ref{comm}.

\section{Chaotic expansions}\label{Swiener}

Before briefly introducing the notion of Wiener chaos, let us recall an integral representation for the nodal length (here we state it for the spherical case but it is -- essentially -- valid in general, thanks to the coarea formula \cite[p.169]{AT}). 

The nodal length on the sphere can be (formally) written as 
\begin{equation}\label{rapSphere}
\mathcal L_\ell = \int_{\mathbb S^2} \delta_0(T_\ell(x)) \|\nabla T_\ell(x)\|\,dx
\end{equation}
where $\delta_0$ denotes the Dirac mass in $0$, $\nabla T_\ell$ the gradient field and $\| \cdot \|$ the Euclidean norm in $\R^2$. 
%
Indeed, let us consider the $\eps$-approximating random variable
$$
\mathcal L_\ell^\eps := \frac{1}{2\eps} \int_{\mathbb S^2} 1_{[-\eps, \eps]}(T_\ell(x))\, \| \nabla T_\ell(x)\|\, dx,
$$
where $1_{[-\eps, \eps]}$ denotes the indicator function of the interval $[-\eps, \eps]$. It is possible to prove that 
\begin{equation*}
\lim_{\eps \to 0} \mathcal L_\ell^\eps = \text{length}(T_\ell^{-1}(0)),
\end{equation*}
both a.s. and in $L^2(\P)$, see \cite{MRW17}, thus justifying \paref{rapSphere}. 
In particular, $\mathcal L_\ell$ is a square integrable functional of a Gaussian field; this is the key point for the theory of chaotic expansions to apply. 

Analogously, for the nodal length of arithmetic random waves we have (with obvious notation)
\begin{equation}\label{rapTorus}
\mathcal L_n = \int_{\mathbb T^2} \delta_0(T_n(x)) \|\nabla T_n(x)\|\,dx. 
\end{equation}
See \cite{MPRW} for details. 

\subsection{Wiener chaos} 

In this part we introduce the notion of Wiener chaos restricting ourselves to our specific spherical setting (the toral case is analogous). For a complete discussion see \cite[\S 2.2]{NP12} and the references therein. 

Denote by $\{H_k\}_{k\ge 0}$ the sequence of Hermite polynomials on $\mathbb{R}$; these polynomials are defined recursively as follows: $H_0 \equiv 1$ and
 $$H_{k}(t) = tH_{k-1}(t) - H'_{k-1}(t), \qquad k\ge 1.$$
Recall that $\mathbb{H} := \{(k!)^{-1/2} H_k, k\ge 0\}$ constitutes a complete orthonormal system in the space of square integrable real functions $L^2(\gamma)$ w.r.t. the standard Gaussian density $\gamma$ on the real line.

Random spherical harmonics \eqref{defSphere} are a by-product of a family of complex-valued Gaussian random variables $\{a_{\ell,m} : \ell = 0, 1, 2, \dots, m=-\ell, \dots, \ell\}$ such that {\bf (a)} every $a_{\ell,m}$ has the form $x_{\ell,m}+iy_{\ell,m}$, where $x_{\ell,m}$ and $y_{\ell,m}$ are two independent real-valued Gaussian random variables with mean zero and variance $1/2$; {\bf (b)} $a_{\ell,m}$ and $a_{\ell',m'}$ are  independent whenever $\ell \ne \ell'$ or $m'\notin \{m,-m\}$, and {\bf (c)} $\overline{a_{\ell,m}} =(-1)^\ell a_{\ell, -m}$. Define the space ${\bf A}$ to be the closure in $L^2(\mathbb{P})$ of all real finite linear combinations of random variables $\xi$ of the form $$\xi = z \, a_{\ell,m} + \overline{z} \, (-1)^\ell a_{\ell, -m},\qquad z\in \mathbb{C},$$ thus ${\bf A}$ is a real centered Gaussian Hilbert subspace
of $L^2(\mathbb{P})$.

{\rm Let us fix now an integer $q\ge 0$; the $q$-th Wiener chaos $C_q$ associated with ${\bf A}$ is defined as the closure in $L^2(\mathbb{P})$ of all real finite linear combinations of random variables of the type
$$
H_{p_1}(\xi_1)\cdot H_{p_2}(\xi_2)\cdots H_{p_k}(\xi_k)
$$
for $k\ge 1$, where the integers $p_1,...,p_k \geq 0$ satisfy $p_1+\cdots+p_k = q$, and $(\xi_1,...,\xi_k)$ is a standard real Gaussian vector extracted
from ${\bf A}$ (in particular, $C_0 = \mathbb{R}$).}

Taking into account the orthonormality and completeness of $\mathbb{H}$ in $L^2(\gamma)$, together with a standard monotone class argument (see e.g. \cite[Theorem 2.2.4]{NP12}), it is possible to prove that $C_q \,\bot\, C_m$ in  $L^2(\mathbb{P})$ for every $q\neq m$, and moreover
\begin{equation*}
L^2(\Omega, \sigma({\bf A}), \mathbb{P}) = \bigoplus_{q=0}^\infty C_q,
\end{equation*}
that is, every real-valued functional $F$ of ${\bf A}$ can be (uniquely) represented as a series, converging in $L^2$, of the form
\begin{equation}\label{e:chaos2}
F = \sum_{q=0}^\infty F[q],
\end{equation}
where $F[q]:={\rm proj}(F \, | \, C_q)$ stands for the the projection of $F$ onto $C_q$ ($F[0]={\rm proj}(F \, | \, C_0) = \E [F]$).

\smallskip


\subsubsection{Random nodal lengths}

Consider now the integral representation \paref{rapSphere} for the nodal length of random spherical harmonics. It can be equivalently written as 
\begin{equation}\label{rapSphere2}
\mathcal L_\ell = \int_{\mathbb S^2} \delta_0(T_\ell(x)) \|\nabla T_\ell(x)\|\,dx
=\sqrt{\frac{\ell(\ell+1)}{2}}\int_{\mathbb S^2} \delta_0(T_\ell(x)) \|\widetilde \nabla T_\ell(x)\|\,dx,
\end{equation}
where  $\widetilde \nabla$ is the normalized gradient, i.e. $\widetilde \nabla:= \nabla / \sqrt{\frac{2}{\ell(\ell +1)}}$ (see \S 3.2.1 in \cite{MRW17} for details). 

We are going to recall the chaotic expansion \paref{e:chaos2} for $\mathcal L_\ell$ 
\begin{equation}\label{chaos_decomp}
\mathcal L_\ell = \sum_{q=0}^{+\infty} \mathcal L_\ell[2q],
\end{equation}
where $\mathcal L_\ell[2q]$ denotes the orthogonal projection of $\mathcal L_\ell$ onto  $C_{2q}$. Note that projections on odd chaoses vanish since the integrand functions in \paref{rapSphere2} are both even.

In \cite[$\S2$]{MRW17} the terms of the series on the r.h.s. of \paref{chaos_decomp} 
are explicitly given (see also \cite{Ros15, MPRW}). 
As it will be clear later in this section, it suffices to deal with the first three terms corresponding to $q=0$ (the mean of the random variable) and $q=1, 2$. See \cite[$\S2$]{MRW17} for more details.

Let us introduce now two sequences of real numbers $\lbrace \beta_{2k}\rbrace_{k=0}^{+\infty}$ and $\lbrace \alpha_{2n,2m}\rbrace_{n,m=0}^{+\infty}$ corresponding to the chaotic coefficients of the Dirac mass at $0$ and the Euclidean norm respectively: for $k=0, 1, 2, \dots$
\begin{equation*}
\beta_{2k} := \frac{1}{\sqrt{2\pi}} H_{2k}(0),
\end{equation*}
while for $n,m=0, 1, 2, \dots$ 
\begin{equation*}
\alpha _{2n,2m}:=\sqrt{\frac{\pi }{2}}\frac{(2n)!(2m)!}{n!m!}\frac{1}{2^{n+m}%
}p_{n+m}\left(\frac{1}{4}\right),
\end{equation*}%
where $p_{N}$ is the swinging factorial coefficient
\begin{equation*}
p_{N}(x):=\sum_{j=0}^{N}(-1)^{j}(-1)^{N}\left(
\begin{array}{c}
N \\
j%
\end{array}%
\right) \frac{(2j+1)!}{(j!)^{2}}x^{j}.
\end{equation*}
The first few terms are
\begin{equation}\label{first_coeff}
\begin{split}
\beta _{0}&=\frac{1}{\sqrt{2\pi }},\quad \text{ }\beta _{2}=-\frac{1}{\sqrt{%
2\pi }},\quad \beta _{4}=\frac{3}{\sqrt{2\pi }},\\
\alpha _{00}&=\sqrt{\frac{\pi }{2}},\quad  \alpha _{02}=\frac{1}{2}\sqrt{%
\frac{\pi }{2}},\quad \alpha _{04}=-\frac{3}{8}\sqrt{\frac{\pi }{2}}.
\end{split}
\end{equation}%
The chaotic expansion of the nodal length is 
\begin{equation}\label{chaos_decomp_sphere}
\begin{split}
\mathcal{L}_{\ell } =\sum_{q=0}^{+\infty} \mathcal L_\ell[2q] = &\sqrt{\frac{\ell (\ell
+1)}{2}}  \sum_{q=0}^{\infty }\sum_{u=0}^{q}\sum_{k=0}^{u}\frac{\alpha
_{2k,2u-2k}\beta _{2q-2u}}{(2k)!(2u-2k)!(2q-2u)!}\times \\
&\times \int_{\mathbb{S}^{2}}H_{2q-2u}(T_{\ell }(x))H_{2k}(\widetilde \partial
_{1;x}T_{\ell }(x))H_{2u-2k}(\widetilde \partial
_{2;x}T_{\ell }(x))\,dx,
\end{split}
\end{equation}
where we use spherical coordinates (colatitude $\theta ,$ longitude $%
\varphi $) and for $x=(\theta _{x},\varphi _{x})$ we are using the notation
\begin{equation*}
\widetilde \partial _{1;x}=\left. (\ell(\ell+1)/2)^{-1/2}\cdot \frac{\partial }{\partial \theta }\right\vert
_{\theta =\theta _{x}},\quad \widetilde \partial _{2;x}= (\ell(\ell+1)/2)^{-1/2}\cdot\left. \frac{1}{\sin \theta }%
\frac{\partial }{\partial \varphi }\right\vert _{\theta =\theta _{x},\varphi
=\varphi _{x}}.
\end{equation*}
It is obvious that analogous formulas as \paref{chaos_decomp_sphere}  hold for the chaotic components of the nodal length $\mathcal L_n$ of arithmetic random waves. 
We are now in a position to give a sketch of the proof of Theorem \ref{thSphere} (see \cite[\S 3]{MRW17} for details).

\section{On the proof of Theorem \ref{thSphere}}\label{Sgreen}

\begin{proof} Substituting \paref{first_coeff} into \paref{chaos_decomp_sphere} we have
\begin{equation}\label{calcoli_media}
\mathcal L_\ell[0] = 4\pi \cdot \frac{1}{2\sqrt 2} \sqrt{\ell(\ell+1)} = \E[\mathcal L_\ell].
\end{equation}
For the second chaotic projection, recalling that $H_2(t) = t^2-1$ and noting that $\alpha_{20}=\alpha_{02}$, we can write 
\begin{equation}\label{2chaos1}
\begin{split}
\mathcal L_\ell[2] = \sqrt{\frac{\ell(\ell+1)}{2}}\Big ( & \frac{\beta_2 \alpha_{00}}{2!}\int_{\mathbb S^2} (T_\ell(x)^2 - 1)\,dx \\
&+ \frac{\beta_0 \alpha_{20}}{2!}\int_{\mathbb S^2} \Big ( \langle \widetilde \nabla T_\ell(x)), \widetilde \nabla T_\ell(x)\rangle - 2\Big )\,dx \Big ).
\end{split}
\end{equation}
A standard application of Green's identity on manifolds\footnote{The author thanks Domenico Marinucci for sharing with her this insight} (see also the proof of \cite[Proposition 7.3.1.]{Ros15}) yields 
\begin{equation}\label{2chaos1}
\begin{split}
\mathcal L_\ell[2] = \sqrt{\frac{\ell(\ell+1)}{2}}\Big ( & \frac{\beta_2 \alpha_{00}}{2!}\int_{\mathbb S^2} (T_\ell(x)^2 - 1)\,dx \\
&+ \frac{\beta_0 \alpha_{20}}{2!}\int_{\mathbb S^2} \Big ( \langle \widetilde \nabla T_\ell(x)), \widetilde \nabla T_\ell(x)\rangle - 2\Big )\,dx \Big )\\
=\sqrt{\frac{\ell(\ell+1)}{2}}\Big ( & \frac{\beta_2 \alpha_{00}}{2!}\int_{\mathbb S^2} (T_\ell(x)^2 - 1)\,dx \\
& +\frac{\beta_0 \alpha_{20}}{2!}\int_{\mathbb S^2} \Big ( - \frac{2}{\ell(\ell+1)}T_\ell(x) \Delta T_\ell(x) - 2\Big )\,dx \Big )\\
=\sqrt{\frac{\ell(\ell+1)}{2}}\Big ( & \frac{\beta_2 \alpha_{00}}{2!}\int_{\mathbb S^2} (T_\ell(x)^2 - 1)\,dx 
+ \frac{\beta_0 \alpha_{20}}{2!}\int_{\mathbb S^2} \Big ( 2 T_\ell(x)^2  - 2\Big )\,dx \Big ),
\end{split}
\end{equation}
where the last equality we used the fact that $T_\ell$ is an eigenfunction of the Laplacian with eigenvalue $-\ell(\ell+1)$. From \paref{2chaos1} we have 
\begin{equation}\label{2chaos2}
\begin{split}
\mathcal L_\ell[2] 
= \sqrt{\frac{\ell(\ell+1)}{2}}\Big (& \frac{\beta_2 \alpha_{00}}{2!}+  \beta_0 \alpha_{20}\Big )\int_{\mathbb S^2} (T_\ell(x)^2 - 1)\,dx = 0,
\end{split}
\end{equation}
where we used \paref{first_coeff}. 
Let us now investigate the fourth chaotic component. To this aim, consider the term 
\begin{equation}\label{Iell}
\mathcal I_\ell:=- \frac{1}{4} \sqrt{\frac{\ell(\ell+1)}{2}}\frac{1}{4!}\int_{\mathbb S^2} H_4(T_\ell(x))\,dx
\end{equation}
whose mean is zero, and whose variance \cite{MR15, MW14} is 
$$
\Var(\mathcal I_\ell) = \frac{1}{16} \frac{\ell(\ell+1)}{2} \frac{1}{4!} 2 \cdot 4\pi \cdot 2\pi \int_{0}^{\pi/2} P_\ell(\cos \theta)^4 \sin \theta\,d\theta.
$$ 
A careful analysis of the fourth moment  of Legendre polynomials (by means of Hilb's asymptotic formula \paref{Hilb}) give, as $\ell\to +\infty$, 
\begin{equation}\label{ciaociao}
\Var \left (\mathcal I_\ell \right ) \sim \frac{1}{32}\log \ell 
\end{equation}
i.e., the variance of this single term is asymptotically equivalent to the total variance \paref{varSphere}.
Investigating asymptotic moments of product of powers of Legendre polynomials and their derivatives we prove that, as $\ell\to + \infty$, 
\begin{equation}\label{opiccolo2}
\frac{\mathcal L_\ell[4]}{\sqrt{\Var(\mathcal L_\ell[4]) }} \text{\quad and\quad} \frac{\mathcal I_\ell}{\sqrt{\Var(\mathcal I_\ell)}}
\end{equation}
are fully correlated. 
By \paref{chaos_decomp}, the orthogonality of Wiener chaoses, \paref{opiccolo2}, \paref{ciaociao}  and \paref{varSphere} we have 
\begin{equation}\label{opiccolo}
\frac{\mathcal L_\ell - \E[\mathcal L_\ell]}{\Var(\mathcal L_\ell)} = \frac{\mathcal I_\ell}{\sqrt{\Var(\mathcal I_\ell)}} + o_{\mathbb P}(1),
\end{equation}
$o_{\mathbb P}(1)$ denoting a sequence converging to zero in probability, 
meaning that the fourth order chaotic component or more precisely the term $\mathcal I_\ell$ ``dominates" the whole series. 
Thus in order to understand the asymptotic behavior of the total nodal length it suffices to study the second order fluctuations of $\mathcal I_\ell$. 
Proposition 3.4 in \cite{MW14} states that, as $\ell\to +\infty$, 
\begin{equation*}
\frac{\mathcal I_\ell}{\sqrt{\Var(\mathcal I_\ell)}}\mathop{\to}^d Z,
\end{equation*}
where $Z$ is a standard Gaussian random variable. This concludes the proof.

\end{proof}

Inspired by \cite{PR}, it should be possible to show directly \paref{opiccolo} by proving in addition that, as $\ell\to +\infty$,  
$$
\Var\left (\sum_{q=3}^{+\infty} \mathcal L_\ell[2q]     \right ) = O(1).
$$

\subsection{Some insight}\label{comm}

\subsubsection{Universality of the mean nodal length}

Consider point (i) in Remark \ref{rem1}. It is immediate to explain this phenomenon in terms of chaotic components, indeed behind \paref{calcoli_media} there is the following formula
$$
\E[\mathcal L_\ell] = \mathcal L_\ell [0] = \sqrt{\frac{\ell(\ell+1)}{2}} \cdot \beta_0 \alpha_{0,0} \cdot \int_{\mathbb S^2} \underbrace{H_0(T_\ell(x))}_{ = 1}\,dx,
$$
where $\ell(\ell+1)/2$ is the variance of the derivatives of $T_\ell$. The same formula, accordingly modified, holds on the torus (and in more general settings -- see also \cite{CM16}). Of course, the same result can be obtained by Kac-Rice formula \cite[Ch. 12]{AT}. 
\subsubsection{Berry's cancellation phenomenon} 
As briefly explained in (ii) Remark \ref{rem1}, Berry's cancellation phenomenon concerns the order of magnitude of the asymptotic variance of the nodal length, which turns out to be smaller than expected. This fact can be expressed in terms of chaotic expansions as the vanishing of the second order chaotic component of the nodal length (see \paref{2chaos2}). Indeed otherwise the variance of the nodal length would have been of order $\ell$, the ``natural" one. 
A careful inspection of \paref{2chaos1} reveals that the steps of the proof of \paref{2chaos2} are independent of the underlying manifold, in particular the same result holds for the toral case (i.e. the second chaotic component of the nodal length for arithmetic random waves vanishes \cite{Ros15, MPRW}). More generally, it is natural to guess that the same cancellation phenomenon appears for the nodal volume on so-called \emph{isotropic} manifolds (compact two-point homogeneous spaces \cite{BM16}) of any dimension (like the hyperspheres \cite{MR15}), and multidimensional tori (see also \cite{Cam17}). 
Moreover, this phenomenon has been observed also for other functionals of nodal sets of Gaussian random eigenfunctions, see \cite{CM16} for a complete discussion. 

\subsubsection{Gaussianity: torus vs. sphere} In \S \ref{Sgreen} we gave a sketch of the proof of Theorem \ref{thSphere}: the second order fluctuations of the spherical nodal length $\mathcal L_\ell$ are Gaussian. 
From its proof it turns out that $\mathcal L_\ell$
and the term $\mathcal I_\ell$ in \paref{Iell} have the same asymptotic behavior. As already said, the asymptotic Gaussianity of 
\begin{equation}\label{magic_term}
\int_{\mathbb S^2} H_4(T_\ell(x))\,dx
\end{equation}
and hence of $\mathcal I_\ell$ was proved in \cite{MW14}. It is now natural to ask for the asymptotic behavior of the corresponding term on the torus, i.e.
\begin{equation}\label{magic_term2}
\int_{\mathbb T^2} H_4(T_n(x))\,dx.
\end{equation}
Indeed, also in the toral case the fourth chaotic component dominates the series of the total nodal length (see \cite{MPRW}). From Lemma 5.2 in \cite{MPRW}  we have 
\begin{equation}\label{magic_term3}
\int_{\mathbb T^2} H_4(T_n(x))\,dx = \frac{6}{\mathcal N_n} \left ( \frac{1}{\sqrt{\mathcal N_n/2}} \sum_{\xi\in \Lambda_n^+} (|a_\xi|^2 - 1) \right )^2 - \frac{3}{\mathcal N_n^2} \sum_{\xi\in \Lambda_n} |a_\xi|^4,
\end{equation}
where $\Lambda_n^+ := \lbrace \xi=(\xi_1, \xi_2)\in \Lambda_n : \xi_2 > 0\rbrace$ if $\sqrt{n}$ is not an integer, otherwise $\Lambda_n^+ := \lbrace \xi=(\xi_1, \xi_2)\in \Lambda_n : \xi_2 > 0\rbrace \cup \lbrace (\sqrt{n},0)\rbrace$. From \paref{magic_term3} the standard CLT applied to 
the following sum of i.i.d. random variables 
$$
\frac{1}{\sqrt{\mathcal N_n/2}} \sum_{\xi\in \Lambda_n^+} (|a_\xi|^2 - 1)
$$
entails that  \paref{magic_term2} is not asymptotically Gaussian. This discussion gives moreover some insight  on the reasons for the nature of the limiting distribution found in Theorem \ref{thTorus}.

\section{Further related work}\label{Sfurther}

\subsection{Monochromatic random waves}\label{RRW}
Consider now the general setting as in \S \ref{BG}. 
For a generic Riemannian metric on a smooth compact surface $\mathcal M$, the eigenvalues $\lambda_j^2$ of the Laplace-Beltrami operator (as in \paref{eq}) are simple and one cannot define a meaningful random model as for the toral or the spherical case, i.e. by endowing   each eigenspace with a Gaussian measure. 
To overcome this, Zelditch introduced  in \cite{Z09} the so-called \emph{monochromatic random waves} as general approximate models of random Gaussian Laplace eigenfunctions defined on manifolds not necessarily having spectral multiplicities, see also \cite{CH16} for more details. The monochromatic random wave on $(\mathcal M,g)$ of parameter $\lambda>0$ is  the random field
\begin{equation}\label{phi}
\phi_\lambda(x) := \frac{1}{\sqrt{\text{dim}(H_{\lambda})}} \sum_{\lambda_j\in [\lambda, \lambda + 1]} a_j f_j(x), \quad x\in \mathcal M,
\end{equation}
where the $a_j$ are i.i.d. standard Gaussian random variables, and
$$
H_{\lambda} := \bigoplus_{\lambda_j \in [\lambda, \lambda + 1]} \text{Ker}(\Delta_g + \lambda_j^2\, \text{Id}),
$$
$\text{Id}$ being the identity operator. The field $\phi_\lambda$ is centered Gaussian and its covariance kernel is 
\begin{equation}\label{covRiem}
K_{\lambda}(x,y) := \Cov(\phi_\lambda(x), \phi_\lambda(y)) = \frac{1}{\text{dim}(H_{\lambda})} \sum_{\lambda_j\in [\lambda, \lambda +1]} f_j(x) f_j(y),\quad x,y\in \mathcal M.
\end{equation}
Now fix a point $x\in \mathcal M$, and consider the tangent plane $T_x\mathcal M$ to the manifold at $x$. The {\it pullback Riemannian random wave} (see \cite{CH16}) associated with $\phi_\lambda$ is the Gaussian random field on $T_x\mathcal M$ defined as
$$
\phi_\lambda^x(u) := \phi_\lambda\left ( \exp_x \left ( \frac{u}{\lambda} \right )\right ),\qquad u\in T_x \mathcal M, 
$$
where $\exp_x : T_x\mathcal M \to \mathcal M$ is the exponential map at $x$. The random field $\phi_\lambda^x$ is  centered Gaussian and from \paref{covRiem} its covariance kernel is 
$$
K_{\lambda}^x(u,v) = K_{\lambda}\left(\exp_x \left ( \frac{u}{\lambda}  \right) , \exp_x \left ( \frac{v}{\lambda}\right ) \right ),\qquad u,v\in T_x \mathcal M.  
$$ 

\begin{definition}[See \cite{CH16}]\label{def scaling}{\rm A point $x \in \mathcal{M}$ is  of {\it isotropic scaling} if, for every positive function $\lambda \mapsto r(\lambda)$ such that $r(\lambda) = o(\lambda)$, as $\lambda\to \infty$, one has that
\begin{equation}\label{limit}
\sup_{u,v \in \mathbb{B}(r(\lambda))  } \left| \partial^\alpha\partial^\beta[  K_{\lambda}^x(u,v) -  2\pi \cdot J_0(\| u-v\|_{g_x})]\right| \to 0,\quad \lambda\to \infty,
\end{equation}
where $\alpha, \beta\in \mathbb{N}^2$ are multi-indices labeling partial derivatives with respect to $u$ and $v$, respectively, $\| \cdot \|_{g_x}$ is the norm on $T_x\mathcal M$ induced by $g$, and $\mathbb{B}(r(\lambda))$ is the corresponding ball of radius $r(\lambda)$ containing the origin. The manifold $\mathcal M$ is of isotropic scaling if every $x\in \mathcal M$ is a point of isotropic scaling, and the convergence in \paref{limit} is uniform on $x$ for each $\alpha, \beta\in \mathbb N^2$. 
}
\end{definition}
It is difficult to prove directly that a point $x$ is of isotropic scaling, except on flat tori, but 
sufficient conditions about the geodesics through $x$ are known (see \cite[\S 2.5]{CH16} and the references therein). The condition
that $\mathcal M$ is a manifold of isotropic scaling is \emph{generic} in the space of Riemannian metrics
on any smooth compact manifold, see \cite[\S 2.5]{CH16} for further details.  Note that one can always choose coordinates around $x$ to have $g_x = \text{Id}$, so that the limiting kernel in \paref{limit} coincides with \paref{cov_berry} up to a factor. 

\subsection{Nodal lengths: recent results}

Let us now fix a \emph{$C^1$-convex body} $\mathcal D\subset \R^2$ (that is: $\mathcal{D}$ is a compact convex set with $C^1$-boundary) such that $0\in \mathring{\mathcal D}$ (i.e., the origin belongs to the interior of $\mathcal D$). 
For every $x\in \mathcal{M}$ we define
$$
\mathcal{Z}_{\lambda,\mathcal D}^x := {\rm length}( (\phi_\lambda^x)^{-1}(0)\cap\mathcal D). 
$$
\begin{theorem}[Special case of Theorem 1 in \cite{CH16}]\label{t:ch}Let $x$ be a point  of isotropic scaling, and assume that coordinates have been chosen around $x$ in such a way that $g_x = {\rm Id.}$ Then, as $\lambda \to \infty$,  
$$
\mathcal{Z}_{\lambda, \mathcal D}^x\mathop{\to}^d \sqrt{2\pi}\cdot {\rm length}(B_1^{-1}(0) \cap \mathcal D),
$$
where $B_1$ is Berry's RWM \paref{cov_berry} with $\lambda_j = 1$. 
\end{theorem}
Theorem \ref{t:ch} states that the local behavior of zeros of monochromatic random waves is universal  (see also \cite{NPR17} for more details). It is not an easy task to find the distribution of the limiting random variable ${\rm length}(B_1^{-1}(0) \cap \mathcal D)$ in Theorem \ref{t:ch}, so let us understand at least what happens when the domain $\mathcal D$ grows to $\R^2$. For $E>0$ let us hence consider
\begin{equation}\label{deflength}
\text{length}(B_{1}^{-1}(0)\cap  \sqrt E \cdot \mathcal D) = \sqrt E \cdot  \text{length}(B_{\sqrt{E}}^{-1}(0)\cap  \mathcal D) =: \sqrt E \cdot \mathcal L_E,
\end{equation}
where $B_{\sqrt{E}}$ is Berry's RWM \paref{cov_berry} with $\lambda_j = \sqrt{E}$. 
 We recall now the main results of \cite{NPR17} concerning the behavior of $\mathcal L_E$, as $E\to +\infty$. Note that \paref{meanlength} and \paref{e:varlength} confirm the results predicted by Berry in \cite{berry2002}. 
\begin{theorem}\label{thlength}
The expectation of the nodal length $\mathcal{L}_E$ is 
\begin{equation}\label{meanlength}
\E[\mathcal L_E] = \text{area}(\mathcal D)\,\frac{1}{2\sqrt{2}}\sqrt{E},
\end{equation}
 the variance of $\mathcal L_E$ verifies the asymptotic relation
\begin{equation}\label{e:varlength}
\Var(\mathcal L_E) \sim \text{area}(\mathcal D)\,\frac{1}{512\pi}\log E, \quad E\to\infty. 
\end{equation}
Moreover, as $E\to \infty$,
\begin{equation*}
\frac{\mathcal L_E - \E[\mathcal L_E]}{\sqrt{\Var(\mathcal L_E)}}\mathop{\longrightarrow}^{d} Z,
\end{equation*}
where $Z\sim \mathscr{N}(0,1)$ is a standard Gaussian random variable. 
\end{theorem}
Note that the mean \paref{meanlength} has the same form as in the toral or spherical case (see Remark \ref{rem1}), and the asymptotic variance \paref{e:varlength} is of lower order than expected -- as predicted by Berry in \cite{berry2002}. Thanks to the symmetry of nodal lines on the two dimensional sphere \cite[\S 1.6.2]{Wig}, the result on the asymptotic variance for the nodal length on $\mathbb S^2$ \paref{varSphere} obtained by Wigman in \cite{Wig} is of course consistent to Berry's prediction (\paref{e:varlength} or \cite{berry2002}).

The proof of Theorem \ref{thlength} relies on chaotic expansions as for the toral and spherical case (see \S \ref{Sgreen}). Here, the second order chaotic component does not vanishes identically but it is possible to prove, by means of Green's identity on manifold, that its contribution is negligible. In view of this, the next chaotic component to study is the fourth one which turns out to dominate the whole series, once proved that the contribution of higher order chaoses is negligible (for details see \cite{NPR17}).

\end{document}